\documentclass[12pt, letterpaper]{article}

\usepackage[utf8]{inputenc}
\usepackage[T1]{fontenc}
\usepackage{mathptmx} 

\usepackage{amsmath}%
\numberwithin{equation}{section}
\usepackage{amsthm}
\usepackage{amsxtra}%
\usepackage{amsfonts}%
\usepackage{amssymb}%

\usepackage{hyperref}
\usepackage{graphicx}
\usepackage{caption}
\usepackage{subcaption}
\hypersetup{colorlinks,breaklinks,
	     linkcolor=blue,urlcolor=blue,
	     anchorcolor=blue,citecolor=blue}
\usepackage{xcolor}

\usepackage[T1]{fontenc}
\usepackage{booktabs}
\usepackage[margin=1in]{geometry}
\usepackage{adjustbox,lipsum}

\usepackage[natbibapa]{apacite} 

\def\degree{^\circ}

\providecommand{\keywords}[1]
{
  \small
  \textbf{\textit{Keywords---}} #1
}

\theoremstyle{plain}

\usepackage{algorithm}
\usepackage{algorithmic}
\usepackage{dsfont}
\newcommand{\x}{{\mathbf x}}

\newcommand{\one}{\mathds{1}}
\renewcommand{\v}{{\mathbf v}}
\renewcommand{\t}{{\mathbf t}}
\renewcommand{\b}{{\mathbf b}}
\newcommand{\p}{{\mathbf p}}
\newcommand{\n}{{\mathbf n}}

\newcommand{\X}{{\mathbf X}}
\newcommand{\N}{{\mathbf N}}
\newcommand{\eps}{\varepsilon}

\usepackage{comment}

\begin{document}

\title{The Virtual Goniometer: A new method for measuring angles on 3D models of fragmentary bone and lithics\thanks{Source code for the virtual goniometer can be found here: \url{https://amaaze.umn.edu/software}}}
\author{Katrina Yezzi-Woodley\thanks{Anthropology, University of Minnesota, \url{yezz0003@umn.edu} (corresponding author)} \and Jeff Calder\thanks{School of Mathematics, University of Minnesota} \and Peter J.~Olver\footnotemark[3]  \and Paige Cody \footnotemark[4] \and Thomas Huffstutler \footnotemark[3] \and Alexander Terwilliger \footnotemark[3] \and Annie Melton\thanks{Anthropology, University of Minnesota} \and Martha Tappen\footnotemark[4] \and Reed Coil\thanks{Sociology and Anthropology, Nazarbayev University, Kazakhstan}\and Gilbert Tostevin \footnotemark[4]}

\date{} 

\maketitle

\begin{abstract}
\begin{small}
The contact goniometer is a commonly used tool in lithic and zooarchaeological analysis, despite suffering from a number of shortcomings due to the physical interaction between the measuring implement, the object being measured, and the individual taking the measurements. However, lacking a simple and efficient alternative, researchers in a variety of fields continue to use the contact goniometer to this day. In this paper, we present a new goniometric method that we call the \emph{virtual goniometer}, which takes angle measurements virtually on a 3D model of an object.  The virtual goniometer allows for rapid data collection, and for the measurement of many angles that cannot be physically accessed by a manual goniometer. We compare the intra-observer variability of the manual and virtual goniometers, and find that the virtual goniometer is far more consistent and reliable. Furthermore, the virtual goniometer allows for precise replication of angle measurements, even among multiple users, which is important for reproducibility of goniometric-based research. The virtual goniometer is available as a plug-in in the open source mesh processing packages Meshlab and Blender, making it easily accessible to researchers exploring the potential for goniometry to improve archaeological methods and address anthropological questions.
\end{small}
\end{abstract}


\section{Introduction}
\label{sec:intro}

Goniometry is an important aspect of archaeological and zooarchaeological analysis. The primary tool for studying angles on objects such as bone fragments or lithics is the pocket, contact goniometer, in essence a metal protractor with a rotating arm.   However, its reliability has come under question and alternative methods have been proposed \citep{dibble1980comparative, morales2015measuring, valletta2020measuring, archer2016still}. Nonetheless, it remains a constituent of current lithic analysis \citep{debert2007raspadita,scerri2016can, muller2016identifying,douglass2018core}, and has expanded into zooarchaeological and taphonomic research to study bone fragmentation. In particular it has been applied to differentiate anthropogenic bone fragmentation from fragmentation by carnivores \citep{alcantara2006determinacion, capaldo1994quantitative, coil2017new, dejuana2011testing, moclan2019classifying, moclan2020identifying, pickering2005contribution, pickering2006experimental}.

The goniometer has a long history.  It was first described by Gemma Frisius (a doctor, mathematician, and cartographer) in 1538 and was derived from the astrolabe, which is the predecessor of the total station, an instrument used for surveying and cartography. The first pocket, contact goniometer was designed in the 1780s to measure the angles on crystals, becoming an essential experimental tool for the burgeoning science of crystallography, enabling accurate measurement of the angles between crystal faces, thereby enabling early researchers to establish the foundations of crystal symmetry classes and generate the classification tables of modern crystallography. Although, at the time, the contact goniometer was revolutionary in its precision, effectively demonstrating the uniformity of crystals, it was largely abandoned in the early 1900s with the advent of yet more accurate methods based on x-ray diffraction \citep{burchard1998history}. 

The contact goniometer first appeared in archaeology when \citet{barnes1939differences} used the instrument to differentiate anthropogenically produced stone tools from naturally-occuring conchoidally fractured rocks. Taphonomic applications began when \Citet{capaldo1994quantitative} used the goniometer to distinguish bones broken by carnivores from those broken by hominins. They extrapolated directly from lithic methods using the goniometer to measure the internal platform angle on models of bone flakes created by taking impressions of notches and associated flake scars on experimentally broken bone. In 2006, \citet{alcantara2006determinacion} introduced a method for identifying actors of breakage by using the goniometer to measure a few fracture angles --- meaning the angle of transition from the periosteal surface to the fracture surface --- on long bone shaft fragments. Prior to this, fracture angles were assessed by eye and categorized as oblique, right, or both, as a means of distinguishing green breaks from dry breaks \citep{villa1991breakage}.
Using the goniometer in the analysis of fragmentary faunal assemblages is gaining traction because it permits researchers to collect seemingly more reliable, quantitative data, and opens the possibility for other avenues of analysis to address questions related to hominin and carnivore interactions at important paleoanthropological sites \citep{coil2017new, dejuana2011testing, dominguez2006new, moclan2019classifying, moclan2020identifying, pickering2005contribution, pickering2006experimental}. 

Although goniometry is claimed to be useful in anthropology, the overall reliability of measurements, and hence subsequent inferences, depends on their accuracy, precision, and replicability. For the contact goniometer to be accurate and precise, the target object must be a size that is compatible with the measuring implement. The contact goniometer was originally designed to measure the angle between intersecting flat surfaces, such as are found on crystals, and is less well adapted to surfaces with curvature and other features, such as cylindrical long bones, or uneven surfaces that are found on bones and lithic flakes. The positioning of the goniometer on the object depends on the user, who must ensure that each arm of the instrument lies flat against and perpendicular to the faces being measured.  Curvature variations can make this placement challenging. Furthermore, sometimes the goniometer cannot access the necessary location for taking the measurement. It might be blocked by adhering matrix or some other surface feature such as a protuberance on a bone fragment or the break on the opposite side of the fragment. Measurements taken by the goniometer can be expected to vary $\pm 5^{\circ}$ \citep{capaldo1994quantitative, draper2011comparison} and extremely sharp edges cannot be measured by a goniometer \citep{dibble1980comparative}, which, in some cases, has a lower limit of $20\degree$. As a result of all these constraints, measurements taken at the edge angle of a single flake or bone fragment will be inconsistent \citep{dibble1980comparative, Johnson2019Testing}, which calls into question the accuracy of any comparative studies.  

Another disadvantage of the contact goniometer is the amount of time required to take each manual angle measurement. As a result, its application to assemblages that number in the thousands can become overwhemlingly time-consuming if not simply too arduous to complete. Thus time and physical constraints may make it impossible to capture sufficiently representative data for detailed and rigorous studies. 

Given the many modern digital tools available, a logical next step would be to turn to automatic, or semi-automatic, computerized means for assessing the relevant angles of scanned bone fragments and stone tools. Recently, in the field of lithic analysis, efforts have shifted in a digital direction \citep{archer2016still, valletta2020measuring, weiss2020lichtenberg, weiss2018variability}. \citet{archer2015diachronic,archer2016still} developed an R package that calculates edge angles based on the thickness of the object at a fixed distance from the edge using basic trigonometry. Their method applies PCA to find the principal axes of the lithic, the first determining its long axis, the second its width, defined as the furthest extent of the object in that direction, and the third its thickness, defined as the distance between corresponding points on each biface. The angle at an edge point is then calculated using the isosceles triangle in the plane perpendicular to the edge, whose apex is the edge point and whose base equals the thickness at a specified distance from the edge point. As pointed out in the code description  \citet{Pop2019Lithics3D}, the function will not work if the plane intersects more than one edge. Therefore, the function depends on the object being of a particular shape, is sensitive to the location of the vertices where measurements are taken, and cannot be easily applied to other tool types or bone fragments which have break edges that are not found on lithics, such as spiral breaks. Because basic trigonometry is used, the angle calculation relies on only three points and is highly sensitive to small topographical changes on the object. The virtual goniometer, however, works on completely general digitized solid objects and uses all the vertices within the patch to define the angle of intersection between the two faces and is therefore not sensitive to small topographical deviations. \citet{valletta2020measuring} developed a mean angle measurement procedure available in a stand-alone software program. This program uses 3D models and creates a mean angle value measurement from a number of selected points along an edge \citep{valletta2020measuring}. They use a cylindrical area that encompasses the entire length of the ridge, or break edge, and average the data along that length to fit two planes on either side of the ridge. This is useful when the edges are straight and uniform. However, many lithic elements and bone fragments do not conform to this ideal. To capture changes along the edge requires a different approach. 

Here we introduce a new method for taking angle measurements, called the \emph{virtual goniometer}. This is a virtual method and can be used to measure angles on 3D models of any object. The virtual goniometer offers a precise and replicable way to measure angles that addresses all of the challenges highlighted above. We provide examples of using the virtual goniometer to measure angles on crystals and stone tools. We assess and compare intra-observer variability in fracture angle measurements taken on bone fragments using a contact goniometer versus the virtual goniometer, and find that the virtual goniometer is far more consistent, performing at least 3.6 degrees ($\pm$ 0.8, 95\% C.I.) more consistently than the contact goniometer. We have also implemented the virtual goniometer as a plug-in in the open source mesh processing software packages Meshlab and Blender, making the method widely available to researchers in the field. 

We will describe the mathematical algorithms used to design the virtual goniometer in Section 2.  Section 3 summarizes the methods used in our particular application, and the materials we have applied it to.  Section 4 presents our results and conducts a statistical analysis of the reliability of three methods of angle measurement in our sample --- manual and two different virtual goniometer protocols.  Perhaps not surprisingly, the virtual goniometer outperforms the manual contact goniometer; moreover, the second virtual process, which we call the ``xyz method'', has remarkably accurate and reproducible measurements in comparison with both of the others.  The final two sections discuss the results and conclusions that can be drawn from our study.

\section{The Virtual Goniometer}

In this section we describe the mathematical algorithms used to design the virtual goniometer.  Its implementation and applications to crystals, lithics, and bones will be described in the following section.

The starting point is a mesh that represents (part of) the bounding surface of a solid object --- a bone fragment, a lithic, a crystal, etc.  
Thus the mesh consists of a reasonably dense sample of points on the surface of the scanned object.  With each mesh point, we also require a unit (length one) vector that points in the outward normal direction to the surface at the point.  The unit outward normals can be constructed from a triangulated mesh by averaging the normals over nearby triangles; for a point cloud one can use local Principal Component Analysis (PCA) or other convenient methods to compute the normal.  Once we have the mesh points and their associated unit outward normals, we do not require any further information such as mesh connectivity to effect our algorithm.

To take an angle measurement at a specified location on the surface, the user must set the location, either by clicking on the mesh representation in the chosen software platform, or by inputting its $(x,y,z)$ coordinates directly.  Both methods are supported in our implementation.  The location is assumed to be on, or at least close to a break edge where an angle measurement has meaning.  Choosing a location far away from the edge can lead to spurious angle measurements that are of no use to the studies for which the virtual goniometer is designed.  The user then specifies a \emph{patch} on the surface, that is, a subset of the vertices, centered at the specified location where the virtual goniometer angle measurement is to be computed.  By a patch of radius $r > 0$, we mean the set of all mesh points that are within distance $r$ of the specified location. In Meshlab, distance means \emph{geodesic distance} measured (approximately) along the surface, and is computed using built in code that is based on Dijkstra's algorithm on a $k$-nearest neighbor graph constructed over the point cloud of vertices.  In Blender, we simply use the Euclidean distance between vertices in the ambient three-dimensional space.  The user specifies the radius $r$ by first clicking on the location and dragging the cursor using the graphical interface to see the resulting patch.  Alternatively, the radius can be a priori specified or entered manually.  

The virtual goniometer algorithm has two main steps. In the first step, we use unsupervised machine learning algorithms, in particular, data clustering techniques, to separate the patch into two regions, each corresponding to one side of the break edge. This part of the algorithm relies on a parameter $\lambda\geq 0$ that is used to ``tune'' the segmentation, when needed. The second step of the algorithm uses Principal Component Analysis (PCA) to find the planes of best fit in each region. The angle is then computed as the angle between the normal vectors of these two planes, and returned to the user along with the mean squared error of the fitting.

To describe the algorithm mathematically, we first note that all vectors are considered as column vectors. There are two key matrices that serve as the input to the algorithm: $\X$ will represent all the vertices in the patch, while  $\N$ will represent their corresponding outward unit normal vectors, while $\lambda $ is the aforementioned tuning parameter.  
 The Virtual Goniometer is summarized in Algorithm \ref{alg:VG}, whose output is the angle measurement $\theta$ in degrees, and the goodness of fit $\eps>0$. 

More precisely, let $n$ be the number of vertices in the patch.  We use $\one=[1,\dots,1]^T$ to denote the column vector with $n$ entries all equal to one, and $\|\x\|$ to denote the Euclidean norm of the vector $\x$. Let $\X$ be the $3\times n$ matrix whose $i^{\rm th}$ column is the vector of $x,y,z$ coordinates of the $i^{\rm th}$ vertex in the selected patch. Let $\N$ be the $3\times n$ matrix  whose $i^{\rm th}$ column is the outward unit normal vector to the surface at the $i^{\rm th}$ vertex.

\begin{algorithm}[tb]
\caption{\ \ Virtual Goniometer}\label{alg:VG}
\begin{algorithmic}[1]
\STATE {\bfseries Input:} Points $\X$, normals $\N$, and tuning parameter $\lambda\geq 0$ 
\STATE {\bfseries Output:} Angle $\theta$ and goodness of fit $\eps$
\STATE $\overline{\x} \gets $ Centroid of vertices $\X$. \ \COMMENT{See notes.}
\STATE $r \gets $ Radius of patch $\X$. \ \COMMENT{See notes.}
\STATE $\t \gets $ Eigenvector of $\N\N^T$ with smallest eigenvalue. \ \COMMENT{Tangent to break curve.}
\STATE $\n \gets\displaystyle \frac{1}{n} \N\one$ \ \COMMENT{Average outward normal vector.}
\STATE $\b \gets \displaystyle \frac{\t\times \n}{\|\t\times \n\|}$ \ \COMMENT{Cross product, normalized --- binormal vector.}
\STATE $\p \gets \b^T\left[\,\N + (\lambda/r)(\X-\overline{\x})\,\right]$ \ \COMMENT{Dot product of all rows with $\b$.}
\STATE $s \gets \text{WithinSS}(\p)$  \ \COMMENT{Find best clustering of 1D data $\p$;  see notes.}
\STATE $\X_- = \X\,[\,\p < s\,], \quad \X_+ = \X\,[\,\p \geq s\,]$   \ \COMMENT{Split into subpatches divided by  break curve.}
\STATE $(\v_\pm,\eps_\pm) \gets \text{PCA}(\X_\pm)$  \COMMENT{Principal component with smallest variance, mean squared error.}
\STATE $\theta \gets 180 - \arccos\left( \text{sign}(\n^T\v_+)\>\text{sign}(\n^T \v_-)\>\v_+^T \v_- \right)$ \ \COMMENT{VG angle --- arccos returns degrees.} 
\STATE $\eps \gets \eps_+ + \eps_-$  \ \COMMENT{Goodness of fit.}
\end{algorithmic}
\end{algorithm}

\begin{figure}[ht] 
    \centering
    \begin{subfigure}[b]{0.15\textwidth} 
        \centering
        \includegraphics[width=\textwidth]{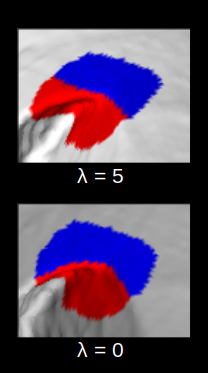}
        \caption{Curve}
        \label{paramA}
    \end{subfigure}  
    \begin{subfigure}[b]{0.148\textwidth} 
        \centering
        \includegraphics[width=\textwidth]{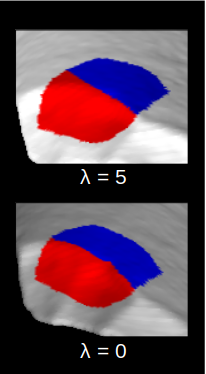}
        \caption{Subtle}
        \label{paramB}
    \end{subfigure} 
    \begin{subfigure}[b]{0.152\textwidth} 
        \centering
        \includegraphics[width=\textwidth]{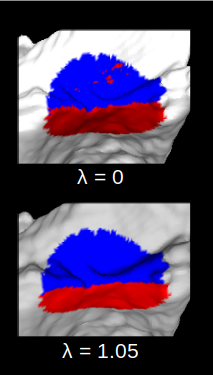}
        \caption{Rugose}
        \label{paramC}
    \end{subfigure} 
        \caption{Examples of how the tuning parameter $\lambda$ affects how the measurement is taken.}   
        \label{paramcomball}
\end{figure}

The Virtual Goniometer Algorithm uses both the normal vectors $\N$ and the vertices $\X$ to segment the patch into two regions. The parameter $\lambda$ controls the tradeoff between how much to rely on the normals versus the vertices. Setting $\lambda=0$ results in clustering using only the normal vectors $\N$. If the surface is noisy, this can give a poor segmentation, since some normal vectors could point in a similar direction even if they are on opposite sides of a break. Increasing $\lambda$ encourages the segmentation to put points that are close together into the same region, and can help to improve the segmentation on noisy meshes. Figure \ref{paramcomball} shows the effect of $\lambda$ on the segmentation, and there are generally three instances, when it  needs to be adjusted: (\ref{paramA}) sharp curves in the ridge; (\ref{paramB}) subtle ridges, usually associated with obtuse angles; and (\ref{paramC}) rugose surfaces. In our implementation, we take  $\lambda=2$ as the default value of the tuning parameter and find that we only need to change this parameter for a small number of measurements.  Moreover, the user can easily learn when and how to adjust it.

Let us describe the individual steps in our algorithm in some detail.
In Step 3, one can use for $\overline{\x}$ any reasonable notion of centroid of the patch $\X$, and either the mean value of the coordinates or a geodesic centroid work well. In our implementation in Meshlab, we set $\overline{\x}$ to be the vertex selected by the user in the click-and-drag selection method. In Step 4, we use the geodesic radius of the patch for $r$; that is the largest geodesic distance from $\overline{\x}$ to any point in $\X$ as measured along the surface. It is also possible to use the Euclidean radius of the patch, and the only difference is a minor change in the effect of the tuning parameter $\lambda$.  However, the algorithm is not overly sensitive to this effect. 

\begin{figure}[t!]
\centering
\includegraphics[width=0.55\textwidth]{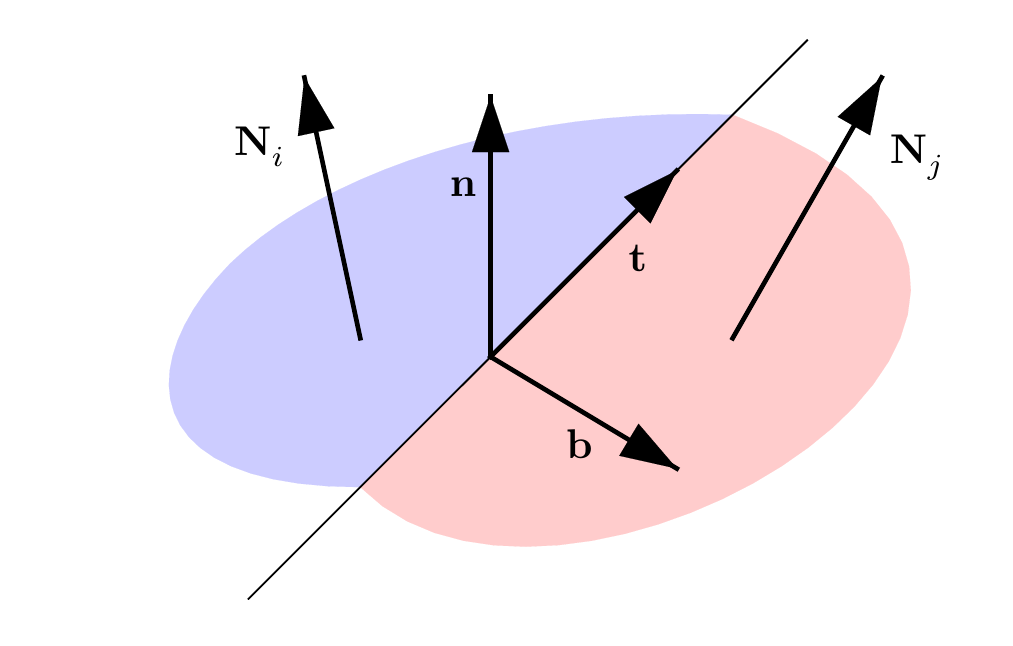}
\caption{\footnotesize{Depiction of a surface patch and the tangent, normal and binormal vectors $\t,\n,\b$ in Algorithm \ref{alg:VG}. The vector $\t$ is tangent to the break curve, which is a line in this example, $\n$ is the average of the outward normals on both sides of the break edge, denoted $\N_i$ and $\N_j$ in the figure, and the binormal $\b$ is orthogonal to $\n$ and $\t$, and hence points \emph{across} the break edge. The basis for our segmentation algorithm is that the sign of the dot product $\b\cdot \N_i$ will indicate on which side of the break edge a particular vertex falls. }}
\label{fig:breakangle}
\end{figure}

The vector $\t$ produced in Step 5 is to be interpreted as the tangent vector to the ``break curve'' that separates the two approximately planar regions on the surface. The vector $\n$ in Step 6 is the averaged outward normal vector to the patch.  Thus, in Step 7, the vector $\b$ given by normalizing the cross product of the break curve tangent and the surface normal is orthogonal to both and hence can be interpreted as the unit \emph{binormal} vector of the patch $\X$, pointing \emph{across} the break curve, i.e., a unit vector that is both tangent to the surface and normal to the break curve.  The binormal vector $\b$ is used to quickly obtain the correct segmentation, as described below. See Figure \ref{fig:breakangle} for a depiction of the vectors $\t,\b$, and $\n$ when the surface is (approximately) planar on either side of the (approximately) straight break edge.

Our clustering method used in Step 9 to divide the surface into two classes was inspired by the random projection clustering methods of \cite{PURR1729,yellamraju2018clusterability}, which involve repeatedly randomly projecting the data to one dimension, and then using the function {\tt WithinSS} described below to perform the clustering of the resulting one-dimensional projected data.  Given one-dimensional data represented by $$\p=(p_1,p_2,\dots,p_n) \qquad \text{where} \qquad  p_1 \leq p_2 \leq \> \cdots \> \leq p_n,$$
this function computes a real number $s$ for which the quantity
\[f(s) = \sum_{p_i\geq s}(p_i-c_1)^2 + \sum_{p_i< s}(p_i-c_2)^2 \qquad \text{where} \qquad c_1 = \frac{\sum_{p_i\geq s}p_i}{\sum_{p_i\geq s} 1},\qquad c_2 = \frac{\sum_{p_i< s}p_i}{\sum_{p_i< s} 1}.\]
is minimized.
The value of $s$ that minimizes $f(s)$ gives the optimal clustering of the one-dimensional data $\p=(p_1,\dots,p_n)$ into two groups $\{p_i \geq s\}$ and $\{p_i< s\}$. Note that $f(s) $ is constant on each interval $p_i < s \leq p_{i+1}$, and hence we can the global minimizer of $f$ simply by computing $f(p_i)$ for $i=1,\dots,n$, and choosing $s = p_i$ that gives the smallest value. (And hence there is also no need to actually sort the data points $\p$.)  This highlights the advantage of working with one-dimensional data; it is very simple and fast to compute optimal clusterings.  

Step 8 in the Virtual Goniometer Algorithm implements a one-dimensional projection of a $\lambda $ weighted combination of the unit normals to the vertices in the patch and the vertices themselves, shifted by the centroid so as to center the patch at the origin. However, the direction of the projection $\b$ is not random but is very carefully chosen as the \emph{binormal} of the patch. The random projection algorithm advocated in \cite{PURR1729} also works very well in this application; however, it requires around $100$ random projections to obtain reliable results and is thus significantly slower. We have also experimented with other clustering algorithms, such as the hyperspace clustering method of \cite{zhang2010randomized}, which also gives very good results at the expense of longer computation times. Our method presented in Algorithm \ref{alg:VG} is very efficient and is suitable for real-time computations in mesh processing software such as Meshlab.

The resulting pair of clusters 
$$\X_- = \X[\p < s] = \{\X_i\mid p_i < s\} \qquad \text{and} \qquad X_+ \X[\p \geq s] := \{\X_i\mid \p_i \geq s\}$$
 contain the vertices belonging to the two components of the surface patch on either side of the break curve.  (However, we do not construct the actual break curve; nor do we require that the centroid or selected location $\overline{\x}$ be thereon.)  In Steps 10 and 11, the PCA function returns each principal component, denoted $\v_\pm$, with the smallest variance and the corresponding eigenvalues $\eps_\pm$, which represents the mean squared error in the fitting. We also experimented with robust versions of PCA (see \cite{lerman2018overview} for an overview), but did not find the results were any more consistent.  Finally, in Step 13, arccos is the inverse cosine function, measured in degrees, not radians, and we use the branch with values between $0$ and $180$ degrees.

\section{Materials and Methods}

In our applications, the process begins by importing a meshed surface into Meshlab. There are two methods used to take virtual goniometer angle measurements at selected points on the surface mesh, which we call the \emph{click-and-drag method} and the \emph{xyz method}.  
For the click-and-drag method, one simply clicks the desired location on the mesh for the measurement (e.g. the break ridge) and then drags to select the area of the surface to be used for the measurement. The chosen area, or surface patch, represents a disk of radius $r > 0$ that is determined by dragging the cursor.   For the xyz method, the user establishes the location of the measurement by inputting the $xyz$-coordinates of the center of the patch. There is also an option to input a specified radius and change the tuning parameter $\lambda $. The plug-in documentation which provides detailed, step-by-step instructions is available at {\tt https://amaaze.umn.edu/}  

\begin{figure}[ht]
    \centering
    \begin{subfigure}[b]{0.305\textwidth} 
        \centering
        \includegraphics[width=\textwidth]{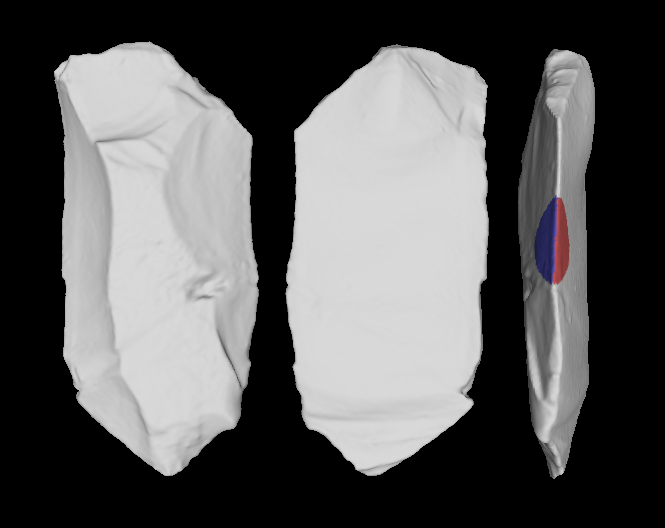}
        \caption{Lithic Flake}
        \label{lithic}
    \end{subfigure}  
    \begin{subfigure}[b]{0.25\textwidth} 
        \centering
        \includegraphics[width=\textwidth]{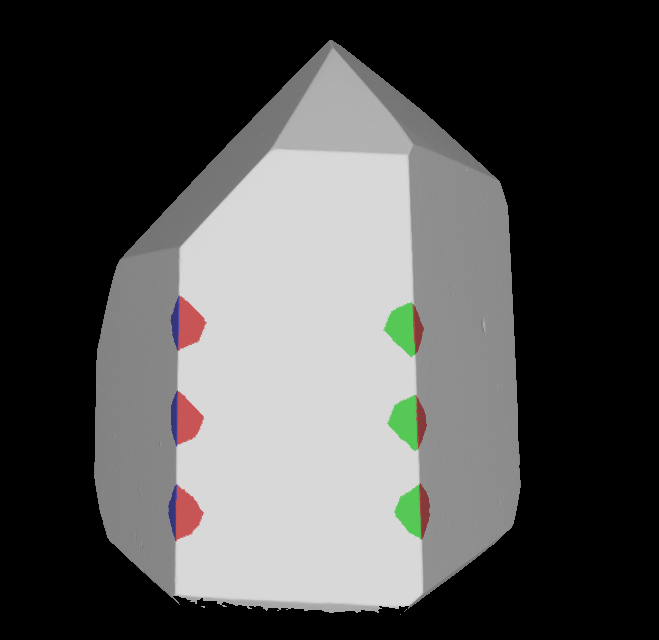}
        \caption{Crystal}
        \label{crystal}
    \end{subfigure} 
    \begin{subfigure}[b]{0.34\textwidth} 
        \centering
        \includegraphics[width=\textwidth]{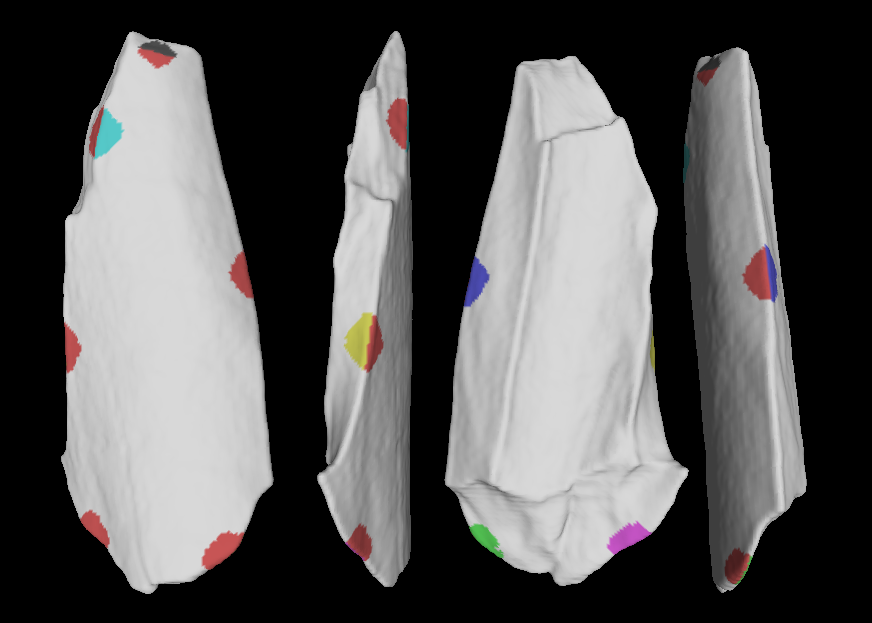}
        \caption{Bone Fragment}
        \label{bone}
    \end{subfigure} 
        \caption{Examples of the Virtual Goniometer on Different Materials}
        \label{examples}
\end{figure}

After the algorithm is run (almost instantaneously), the selected patch is color-coded into two contrasting colors showing the two vertex clusters representing the two sides of the break curve that are used to determine the PCA planes, and hence the angle based on their normals.  The colors rotate through a pallet list of pairs of contrasting colors; see Figure \ref{examples} for examples of the resulting visual output in Meshlab. Further measurements can be taken at or near the original location; see Figure \ref{crystal}. The user can also choose to advance to a new location at which point the colors change; see Figures \ref{crystal} and \ref{bone}. Numerical data for each measurement are automatically recorded in a {\tt .csv} file. 
 
We compared the manual and virtual goniometers for computing fracture angles on a sample of bone fragments consisting of 537 breaks taken from 86 appendicular long bone shaft fragments ($\geq$ 2 cm maximum dimension)  randomly chosen from a collection of experimentally broken \textit{Cervus canadensis nelsoni} (Rocky Mountain elk) and  \textit{Odocoileus virginianus} (white-tailed deer) limb bones. Fragments were scanned using a medical CT scanner at the University of Minnesota's Center for Magenetic Resonance Research (CMRR) (slice thickness: 0.6, reconstruction interval: 0.6 mm, KV: 80, MA:28, rotation time: .05 sec, pitch: 0.8, alogrithm: bone window, convolution kernel: B60f sharp) and then surfaced using MATLAB. The break surfaces on each fragment were manually subdivided into separate break planes \citep{gifford1989ethnographic, haynes1983frequencies, pickering2005contribution}. All breaks were measured, and we did not impose a minimum break length requirement. We measured the fracture angle --- which is defined as the angle of transition between the periosteal surface and break surface --- on each break face of each fragment \citep{alcantara2006determinacion, capaldo1994quantitative, villa1991breakage}. Following the method established by \citet{alcantara2006determinacion} and further described by \citet{pickering2005contribution}, we chose to measure at the center along the break length. 

Each break was measured using both a contact goniometer ("man") and the virtual goniometer using our Meshlab plug-in. Two methods were employed using the virtual goniometer: the click-and-drag method (``drag") and the xyz-coordinate input method (``xyz"). The fragments were first measured using the drag method. After the first round of drag angle measurements, denoted $\theta _{drag}$, screen shots were taken of all the colorized models to create a 2D map of the measurements that served as a visual guide for subsequent drag angle measurements, which we denote by $\varphi _{drag}$ and $\psi _{drag}$.  Three manual angle measurements, $\theta _{man}, \varphi _{man}, \psi _{man}$ were then taken using the contact goniometer, again using the map as a guide.  

The xyz method was used to replicate the first set of drag measurements, so $\theta _{xyz} = \theta _{drag}$. The values for the radius and location (i.e. xyz-coordinates) of $\theta _{drag}$ were input into the plug-in for the xyz method to produce two further angle measurements $\varphi _{xyz}$, $\psi _{xyz}$).  Thus, the xyz method was only executed twice, and, in total, the same person measured the angle of each break eight times (3 manual, 3 drag, and 2 xyz) thereby allowing us to test for intra-observer error. Using additional data automatically provided by the two virtual methods, we analyzed the degree to which the location of the measurement, the distance from the edge, and the number of vertices used impact how the angle measurement varies. 

\section{Results}

\subsection{Summary statistics for all methods}

Of the 537 breaks in our randomly selected sample of bone fragments, 500 (93.1\%) could be measured manually using the contact goniometer, while the other 37 breaks could not be physically measured. For 34 of those breaks, one or both arms of the contact goniometer were blocked from contacting the break face or periosteal face. The 3 remaining breaks that could not be measured manually broke off of a fragment that suffered lab damage after being scanned. The 3D mesh constructed prior to the lab break made it possible to use these measurements for analysis of the virtual goniometer. When comparing the manual measurements to the click-and-drag and xyz methods, we only used the 500 breaks that could be measured by all three methods. 

To test variability, we computed three angle measurements for each of the three methods, which we denote by $\theta, \varphi , \psi$ and use subscripts man, drag, xyz to indicate which method is employed. When using the manual and drag method, the location of the measurement is selected by eye. For the xyz method, the user enters the $x,y,z$ coordinates directly.

To assess how much each method varied, we calculated the Intra-Observer Variability (IOV) for the three angle measurements $(\theta, \varphi, \psi)$ for each break under the three different methods. The IOV is the average of the absolute value of the differences between each of the three angle measurements, all taken at the same location:
\begin{equation}\label{eq:Laplace}
\begin{aligned}
IOV_{man} &= \frac{|\theta _{man}-\varphi _{man}|+|\theta _{man}-\psi _{man}|+|\varphi _{man}-\psi _{man}|}3,\\
IOV_{drag} &= \frac{|\theta _{drag}-\varphi _{drag}|+|\theta _{drag}-\psi _{drag}|+|\varphi _{drag}-\psi _{drag}|}3,\\
IOV_{xyz} &= \frac{|\theta _{xyz}-\varphi _{xyz}|+|\theta _{xyz}-\psi _{xyz}|+|\varphi _{xyz}-\psi _{xyz}|}3.
\end{aligned}
\end{equation}
Ideally, there should be no variation, so that $IOV=0$. 

Table \ref{tab:sumstatsman} shows summary statistics, such as the mean and median values, for the IOV for the three methods. The median for the manual IOV (4.67$^{\circ}$) is marginally better than the expected variation (5$^{\circ}$) described in \citet{capaldo1994quantitative} and \citet{draper2011comparison} but the mean for the manual IOV (7.08$^{\circ}$) is over 2$^{\circ}$ higher than expected and the standard deviation (8.48$^{\circ}$) is high. All but seven of the manual IOVs are $<$31$^{\circ}$. The remaining seven are $>$50$^{\circ}$ and could be considered outliers to which the mean and standard deviation are sensitive. However, removing those seven IOVs would not sufficiently reduce these values because the median IOV for the drag method (2.28$^{\circ}$) and the standard deviation (3.34$^{\circ}$) are considerably smaller. The xyz method has a consistently smaller IOV compared to the other methods, with a median of 0.001$^{\circ}$, mean 0.006$^{\circ}$, and standard deviation 0.011$^{\circ}$.

\vspace{3mm}
\begin{table*}[ht]
\vspace{-3mm}
\caption{Summary statistics for angle IOV$^{\circ}$}
\vspace{-3mm}
\label{tab:sumstatsman}
\vskip 0.15in
\begin{center}
\begin{small}
\begin{sc}
\begin{tabular}{llll}
\toprule
Statistics&Manual&Drag&xyz\\
\midrule
N&500&500&500\\
min&0&0.04&0\\
mean&7.08&3.27&0.006\\
median&4.67&2.28&0.001\\
max&73.33&23.15&0.06\\
sd&8.48&3.34&0.011\\
\bottomrule
\end{tabular}
\end{sc}
\end{small}
\end{center}
\vskip -0.1in
\end{table*}

\begin{figure}[ht]
    \centering
    \begin{subfigure}[b]{0.32\textwidth} 
        \centering
        \includegraphics[clip=true,trim=10 1 30 30, width=\textwidth]{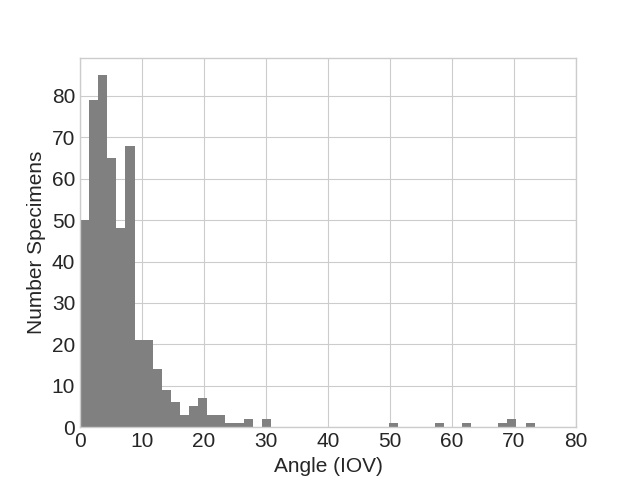}
        \caption{manual}
        \label{histogram_man}
    \end{subfigure}  
    \hfill 
    \begin{subfigure}[b]{0.32\textwidth} 
        \centering
        \includegraphics[clip=true,trim=10 1 34 30, width=\textwidth]{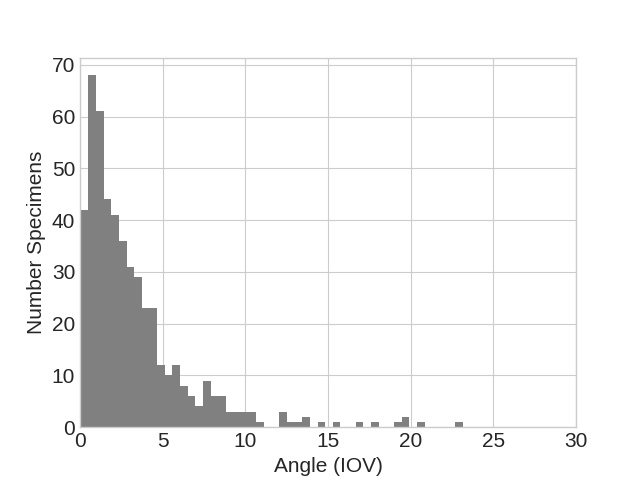}
        \caption{click-and-drag}
        \label{histogram_drag}
    \end{subfigure} 
    \hfill   
    \begin{subfigure}[b]{0.32\textwidth} 
        \centering
        \includegraphics[clip=true,trim=10 1 30 30, width=\textwidth]{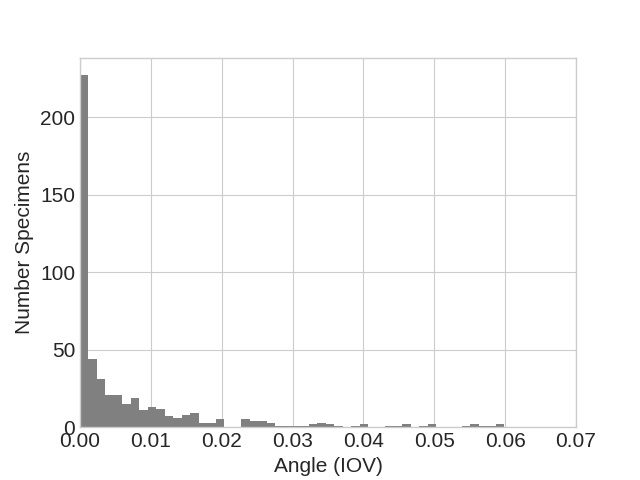}
        \caption{xyz-coordinate}
        \label{histogram_xyz}
    \end{subfigure} 
        \caption{Histograms of angle IOV$^{\circ}$ using all three methods (n=500)}   
        \label{Hist_all}
\end{figure}

Figure \ref{Hist_all} shows histograms of the IOV for each method.  We see that the IOV for the manual and drag method are rather dispersed with a larger proportion of breaks characterized by larger errors, while the IOV for the xyz method is highly concentrated around very low variabilities. We point out to the reader that the scale of the histogram axes are different in each case. The much smaller range and limited dispersion suggests that the virtual goniometer, regardless of method, outperforms the contact goniometer and the xyz method is far more precise than the other two methods.

\subsection{Confidence intervals}

\begin{figure}[ht]
\centering
\includegraphics[width=\textwidth]{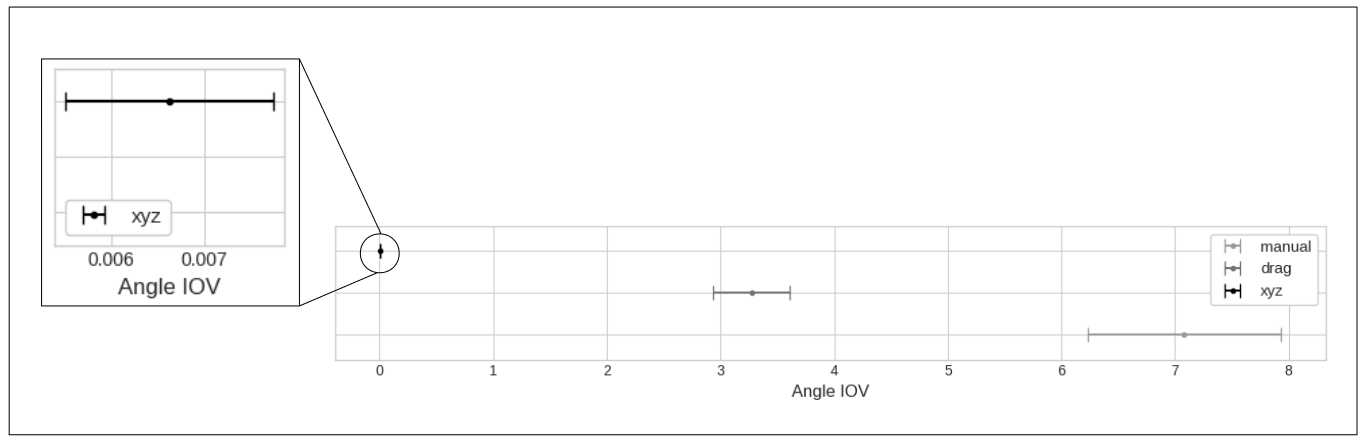}
\caption{Plot of the Confidence Intervals for the IOV's (all methods).}
\label{cipopup}
\centering
\end{figure}

We calculated 95\% and 99\% confidence intervals, the results of which can be found in Table \ref{tab:ci} and Figure \ref{cipopup} \citep{weisberg2005applied}. Even with an increase in the confidence level to 99\%, the difference is striking --- indeed, we have to apply a magnification to observe its confidence interval. The range for the manual method is 2.5 times larger than the drag method and between 760.5 and 772.9 times larger than the xyz method. The drag method is a little over 300 times larger than the xyz method. The small range IOV of the xyz method indicates the method is exceptionally precise, especially compared to the drag and manual methods.

\vspace{3mm}
\begin{table*}[ht]
\vspace{-3mm}
\caption{Confidence Intervals}
\vspace{-3mm}
\label{tab:ci}
\vskip 0.15in
\begin{center}
\begin{small}
\begin{sc}
\begin{adjustbox}{max width=\textwidth}
\begin{tabular}{|ll|lll|lll|}
\hline
&\multicolumn{1}{c}{}&\multicolumn{3}{|c|}{95\% CI}&\multicolumn{3}{c|}{99\% CI}\\
\hline
Method&Mean&Interval&Range&Mean Error&Interval&Range&Mean Error\\
\hline
Manual&7.08&(6.2298, 7.9302)&1.7004&$\pm$0.74345&(6.0153, 8.1447)&2.1294&$\pm$0.62207\\
Drag&3.2740&(2.9388, 3.6093)&0.6705&$\pm$0.29318&(2.8541, 3.6939)&0.8398&$\pm$0.24531\\
xyz&0.0066&(0.0055, 0.0077)&0.0022&$\pm$0.00097&(0.0052, 0.0080)&0.0028&$\pm$0.00082\\
\hline
\end{tabular}
\end{adjustbox}
\end{sc}
\end{small}
\end{center}
\vskip -0.1in
\end{table*}

We ran Tukey's HSD (honestly significant difference) test with confidence level $\alpha=0.05$ to assess the significance of the difference in means of the IOV scores for each method \citep{barnette1998tukey}. Table \ref{tab:ttests} shows that the virtual goniometer's drag method is 3.6 degrees ($\pm$ 0.8, 95\% C.I.) more consistent than the contact goniometer and the virtual goniometer's xyz method is 7.1 degrees ($\pm$ 0.8, 95\% C.I.) more consistent than the contact goniometer. Furthermore, the xyz method is 3.4 degrees ($\pm$ 0.8, 95\% C.I.) more consistent than the drag method. The HSD test rejected the null hypothesis that any of the mean IOV values are equal with statistical significance $p < 0.001$. These methods are not equally effective. It is clear that the xyz method is far superior to both the manual and the drag methods.

\vspace{3mm}
\begin{table*}[ht]
\vspace{-3mm}
\caption{Results of the Tukey HSD Test}
\vspace{-3mm}
\label{tab:ttests}
\vskip 0.15in
\begin{center}
\begin{small}
\begin{sc}
\begin{tabular}{lllllll}
\toprule
group1 & group2 & meandiff & p-adj & lower & upper & reject\\
\midrule
drag IOV & man IOV & 3.637 & 0.001 & 2.873 & 4.402 & True\\
xyz IOV & drag IOV & 3.436 & 0.001 & 2.685 & 4.187 & True\\
xyz IOV  & man IOV & 7.073 & 0.001 & 6.309 & 7.838 & True\\
\bottomrule
\end{tabular}
\end{sc}
\end{small}
\end{center}
\vskip -0.1in
\end{table*}

\subsection{Summary statistics for the IOV (virtual methods)}

\begin{figure}[ht]
    \centering
    \begin{subfigure}[b]{0.45\textwidth} 
        \centering
        \includegraphics[clip=true,trim=10 2 30 30, width=\textwidth]{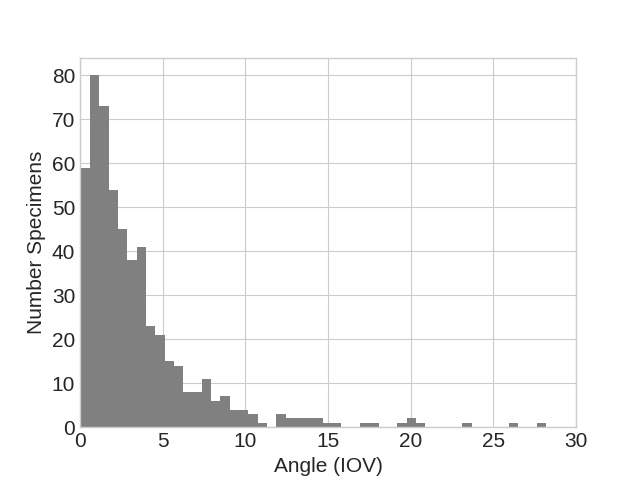}
        \caption{click-and-drag}
        \label{histogram_drag_all}
    \end{subfigure}  
    \hfill 
    \begin{subfigure}[b]{0.45\textwidth} 
        \centering
        \includegraphics[clip=true,trim=10 2 30 30, width=\textwidth]{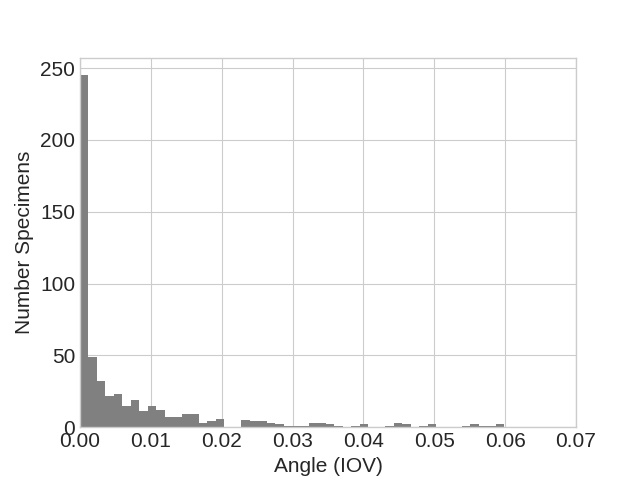}
        \caption{xyz-coordinate}
        \label{histogram_xyz_all}
    \end{subfigure} 
        \caption{Histograms of angle IOV using all virtual methods (n=537)}   
        \label{Hist_virt}
\end{figure}

Since the center, radius, and vertex data cannot be collected using the manual method, we compare the drag and xyz methods for the entire sample using calculations of the IOV for each variable. The IOV of the angle measurement for the drag method features numerous large values and has a standard deviation of 3.7$^{\circ}$, varying as much as 28.2$^{\circ}$ (see Figure \ref{Hist_virt}). The xyz method has a standard deviation of 0.01$^{\circ}$ with a maximum variation of 0.06$^{\circ}$. Most of the xyz IOVs fall below 0.02$^{\circ}$. It is clear that the variation in the angle IOV is the result of variation in the location (represented by the center of the patch), the patch's radius, and the number of vertices in the patch (see Table \ref{tab:vtsummary}).

Table \ref{tab:vtsummary} gives summary statistics for the IOV for the drag and xyz methods based on IOVs for the angle measurement, the number of vertices in the selected patch, the radius of the patch, and the center location of the patch.

\vspace{3mm}
\begin{table*}[ht]
\vspace{-3mm}
\caption{Summary statistics for IOV$^{\circ}$}
\vspace{-3mm}
\label{tab:vtsummary}
\vskip 0.15in
\begin{center}
\begin{small}
\begin{sc}
\begin{tabular}{|l|llll|llll|}
\hline
&\multicolumn{4}{|c|}{Click-and-Drag}&\multicolumn{4}{c|}{xyz-Coordinates}\\
\hline
Statistics&Angle&Vertices&Radius&Center&Angle&Vertices&Radius&Center\\
\hline
N&537&537&537&537&537&537&537&537\\
min&0.04&2.67&0.01&0&0&0&0&0\\
mean&3.44&270.25&0.52&1.74&0.01&0.39&0.00&0\\
median&2.32&170.67&0.41&1.09&0.00&0.67&0.00&0\\
max&28.20&2250&2.63&30.49&0.06&0.67&0.03&0\\
sd&3.70&295.06&0.40&2.41&0.01&0.33&0.00&0\\
\hline
\end{tabular}
\end{sc}
\end{small}
\end{center}
\vskip -0.1in
\end{table*}

\subsection{Multiple regression}

To better understand how changes in the location, the radius, and the number of vertices in the patch affect the angle measurement we ran a multiple regression using a log transformation \citep{weisberg2005applied}.

\begin{figure}[ht]
    \centering
    \begin{subfigure}[b]{0.4\textwidth} 
        \centering
        \includegraphics[clip=true,trim=0 22 30 80, width=\textwidth]{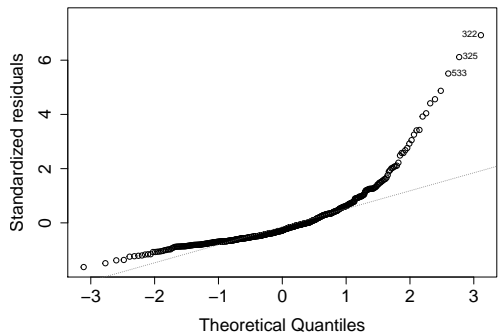}
        \caption{Q-Q plot of the original data}
        \label{QQreg}
    \end{subfigure}  
    \hfill 
    \begin{subfigure}[b]{0.4\textwidth} 
        \centering
        \includegraphics[clip=true,trim=0 22 30 80, width=\textwidth]{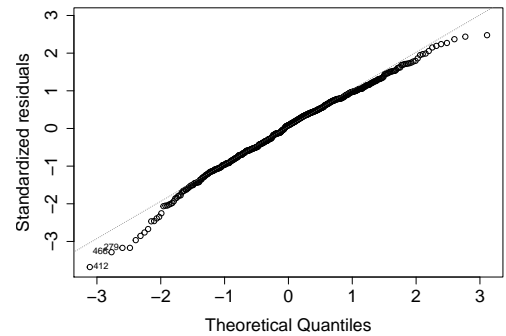}
        \caption{Q-Q plot of the log transformed data}
        \label{QQlog}
    \end{subfigure} 
        \caption{Q-Q plots of the residuals(n=537)}   
        \label{QQ}
\end{figure}

\begin{figure}[ht]
\vspace{4mm}
\centering
\includegraphics[clip=true,trim=0 22 30 80, width=.5\textwidth]{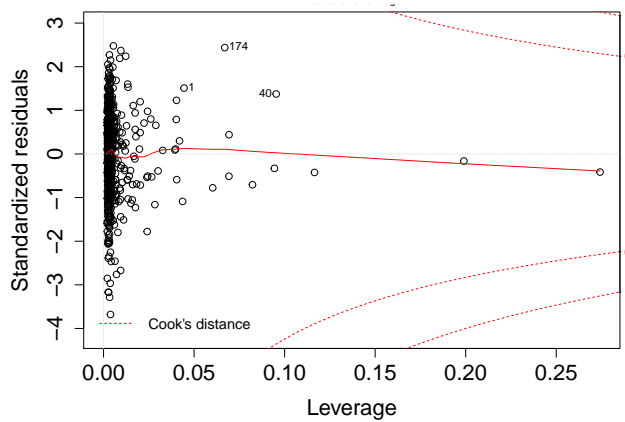}
\caption{Cook's Plot (log-transformed)}
\label{cook}
\centering
\end{figure}

\begin{table}[ht]
\centering
\caption{Multiple Regression Results (n=537)}
\label{mreg}   
\begin{tabular}{rrrrr}
  \hline
& Estimate & Std. Error & t value & Pr($>$$|$t$|$) \\
  \hline
(Intercept) & 0.3995 & 0.0717 & 5.58 & 3.94e-08 \\
  vertices IOV & -0.0009 & 0.0003 & -3.6 & 0.000346 \\
  radius IOV & 0.9997 & 0.1993 & 5.015 & 7.23e-07 \\
  location IOV & 0.0573 & 0.0197 & 2.909 & 0.003772 \\
   \hline
\end{tabular}\\
\vspace{1mm}
\footnotesize{Adjusted R-squared: 0.0872, p-value: 2.242e-10}  
\end{table}

We chose to log-transform the response variable, using the standard natural log, in order to make the data satisfy the assumption of normal error terms (see Figure \ref{QQ}). On the original data the residuals did not follow a normal distribution and were quite skewed (Figure \ref{QQreg}). The residuals for the log-transformed data follow a normal distribution (Figure \ref{QQlog}).  Though the residuals do not have an obvious pattern, the Cook's distance plot shows that none of the datapoints are overly influencing the model (Figure \ref{cook}) \citep{cook1977detection}.  None of the observations have high Cook's values which is indicated by the absence of datapoints in the upper or lower right-hand region of the plot. 

The results of the multiple regression are presented in Table \ref{mreg}. The p-value of the model's fit is significant (p=2.242e-10) as are the p-values for each of the explanatory variables indicating that all three variables are influencing the angle IOV. In this model, the radius (estimate = 0.9997) has the biggest impact. The estimate for the vertices (-0.0009) shows that this has the least impact. The data are highly scattered (R$^2=$0.0872) so it would be difficult to predict the angle IOV based on this model. Nonetheless, it is clear that, when replicating measurements, changes in the location, the radius, and the number of vertices can influence the angle measurement. Having a large number of vertices within a consistently sized patch with a consistent location is the best method for achieving small IOVs which is why the xyz method is the best option for replication of measurements.   

\section{Discussion}
    
The salient questions are whether or not goniometry can be useful for addressing anthropological questions and, if so, in what ways. To answer these questions, measurement methods must result in values that accurately capture the angle of interest with reasonable precision. Methods that do not accurately reflect the desired information undermine the ability to make credible anthropological inferences.

Being able to take a measurement is the first priority. The primary issue with the contact goniometer is the physical interaction between the instrument and its target. When measuring bone fragments it can be difficult to reach the desired location for the measurement. This is often a matter of scale in that the goniometer is too large relative to the size and shape of the bone fragment. While a tiny goniometer would, at least in principle, be able to measure smaller or hard-to-reach locations, in practice it would be difficult for the user to handle comfortably and accurately. In some instances, the arms are too long and are blocked by features on the bone fragment or the break on the opposite side. The latter happens when more than half of the circumference of the fragment is present or the angle of interest is acute and descends into the medullary cavity. When measuring notches, \citet{capaldo1994quantitative} chose to take molds because the goniometer was unable to reach the notch surfaces. However, making molds is often not  an option and can even damage the specimen.   

Even if an acute angle is accessible, some goniometers cannot measure anything less than 20$^{\circ}$ because the arms overlap, blocking the location where the specimen is supposed to fit. To test the accuracy of the virtual goniometer, we created synthetic meshes of intersecting flat planes with known angles. Given a sufficient number of vertices, which depends on scan resolution and radius size, the virtual goniometer can accurately measure an angle as small as 1$^{\circ}$, which is far superior to any version of the contact goniometer in everyday use.  

The only instance in which the virtual goniometer cannot measure a fracture angle is when a sharp natural curve on the bone fragment is close to the fracture edge. This is only an issue when the radius is large. Reducing the radius such that it captures more of the ridge and less of the natural curve rectifies the problem. A smaller radius might be better when measuring angles on bone fragments. Limiting the radius such that it is local to the transition but still sufficiently large to not be hampered by imaging artifacts could provide a more informative measure. As the radius increases in size, it captures changes in the topography. Though the topographical changes on the break surface may also be of interest as it pertains to examining bone fragments, this should remain separate from the angle of transition between faces.   
 
The ability to take a measurement does not ensure its appropriateness. Because the surfaces are not flat, the angle value depends on the distance from the edge where the measurement is taken on both faces. Stopping 3 mm from the edge of an object or stopping 5 mm from the edge can result in a different angle value. When using the contact goniometer, the distance from the edge cannot be controlled. Due to variations in topography among break faces on bone fragments, the arm of the goniometer does not consistently make contact on the surface at the same distance from the edge. When the surface is concave, the goniometer will not rest against the concavity, and will reach across that expanse to the other side. 

Conversely, a convex surface prevents the goniometer from reaching the other side of the break. \Citet{capaldo1994quantitative} bypassed this by cutting into the molds featuring convex surfaces so the goniometer could reach the full extent of the surface, which is not a viable solution when measuring the object directly. In this case, measuring an angle that extends across the whole surface is not possible. When the surface is rugose, the goniometer will make contact with the highest point within its reach at that location. Such topographical variations on the surface of the break prevent the goniometer from connecting with most of the surface and it fails to capture a significant amount of intermediate information at a single location. 

If we were uncertain about a break, we flagged it and categorized it as concave, hinged, surface, edge, or other. If we were unable to measure the break, it was categorized as blocked or other (see Table \ref{tab:category}). Of the flagged fragments, 108 were measured. Arguably, these measurements do not accurately capture the fracture angle and may not be useful for fracture edge analysis.\\ 

\begin{table*}[ht]
\vspace{-3mm}
\caption{Categorized and uncategorized breaks}
\vspace{-3mm}
\label{tab:category}
\vskip 0.15in
\begin{center}
\begin{small}
\begin{sc}
\begin{tabular}{lll}
\toprule
Category&Count&Percent\\
\midrule
All breaks&537&100\%\\
Not categorized&392&73\%\\
Categorized&145&27\%\\
\hspace{5mm}Blocked$^{1}$ &34&6.33\%\\
\hspace{5mm}Concave&70&13.04\%\\
\hspace{5mm}Hinge&21&3.91\%\\
\hspace{5mm}Surface$^{2}$&4&0.74\%\\
\hspace{5mm}Edge$^{3}$&11&2.05\%\\
\hspace{5mm}Other&5&0.93\%\\
Categorized (measured)&108&20.11\%\\
Categorized (not measured)&37&6.89\%\\
\bottomrule
\end{tabular}\\
\tiny{$^{1}$The goniometer arm cannot make contact with the desired surfaces}\\ 
\tiny{$^{2}$The surface is rugose and has a lot topographical relief}\\
\tiny{$^{3}$The edge is rounded or hinders efforts to measure}\\ 
\end{sc}
\end{small}
\end{center}
\vskip -0.1in
\end{table*}

Not only does the surface of a break change across the width of the break, angle values can vary along the length of the break edge. Some researchers measure at the midpoint along the ridge \citep{pickering2005contribution, coil2017new, dibble1980comparative} whereas others use the most extreme angle \citep{moclan2019classifying}. \citet{capaldo1994quantitative} define a midpoint, but, it only applies to notches and cannot be extrapolated to all fracture edges because it depends on features that are specific to notches, specifically inflection points. 

Finding the midpoint on a break edge is less clear. The midpoint on a bone fracture edge as described by \citet{pickering2005contribution} and \citet{coil2017new} is likely the approximate midpoint as opposed to an exact midpoint. The break edge on each break face on a bone fragment can be viewed as a contour. Establishing the endpoints for the contour is the first challenge. Typically, the full length of the break face does not terminate at the same location that the ridge between the periosteal surface and break surface terminates. Once a decision is made as to where the endpoints of the contour will be located, then one can decide where to take the angle measurement. In regard to the midpoint, one could choose the midpoint of the Euclidean distance between the two end points of that contour or one could choose the midpoint of the full contour length. In either case, finding the exact midpoint on the physical object is challenging, if not impossible, and time consuming. For example, calipers could be used in the case of the Euclidean midpoint. However, this would require a consistent orientation of the specimen in relationship to the calipers. Finding the most extreme angle would require approximation as well unless one measures many angles along the edge in order to identify which one is the most extreme. This requires that the angle can be measured at all locations and the goniometer is not blocked by features on the specimen at any location.

When using 3D models, specific points can be chosen as endpoints for the contour. The Euclidean distance or the contour length can be calculated and a midpoint can be consistently defined and extracted. Though this is beyond the purview of the virtual goniometer, the xyz-coordinates for the exact midpoint can be input into the virtual goniometer. 

Whether choosing the most extreme measure or the midpoint, taking a measurement at a single location is anthropologically arbitrary. Bones are not limited to one instance of fragmentation. If the objective is to identify the first actor of breakage, then the extrema or midpoints will not be comparable among specimens if additional fragmentation takes place. An alternative approach, that could provide a richer dataset, would be to take multiple measurements along the edge of the break which can be done quickly using the virtual goniometer. 

Assuming the angle can be taken and is appropriate to take, it also needs to be precise enough that it is useful. When measuring platform angles of notch molds using the contact goniometer, \citet{capaldo1994quantitative} expected measurements to vary up to 5$^{\circ}$. This may not be precise enough. For example, \Citet{alcantara2006determinacion} state that, in general, carnivores produce fracture angles between 80$^{\circ}$ and 110$^{\circ}$, whereas fracture angles on bones broken by percussion will be <80$^{\circ}$ and >110$^{\circ}$, offering a 30$^{\circ}$ range for assigning carnivores as the actors responsible for breakage. If one chooses a more stringent approach where fracture angles between 85$^{\circ}$ and 95$^{\circ}$ \citep{pickering2005contribution}, generally labeled as dry breaks, are excluded, this diminishes the ranges for carnivores to 5$^{\circ}$ and 15$^{\circ}$. These are small windows for such a large error range. In fact, one of those ranges is equal to the expected error range. In a test that looked at the reliability of the contact goniometer when measuring knee angle in a medical context, \citet{draper2011comparison} noted that in order to remain within an error range of 5$^{\circ}$, the location of the measurement had to be within 2mm of the actual center of the patella. Not only is this equal to the error range presented by Capaldo and Blumenschine, these results also highlight the importance of the location where the measurement is taken.  

Being able to control how a measurement is taken is key to extracting anthropologically useful information. The virtual goniometer does not suffer from the limitations of the contact goniometer. It is flexible and can better meet the specific needs of the user who can select the location and which portion of the object to measure. It captures the entire surface in the specified region. The radius (i.e. the distance from the edge) of the patch can be specified. The patch can be selected by eye using the click-and-drag method or by specifying the parameters through numerical inputs. The option to select the location by clicking and then specifying the radius is also available when a consistent radius is required. The segmentation parameter can also be adjusted in order to tweak how the virtual goniometer segments the specified patch into two planes.  

A final consideration is the individual taking the measurement. The way in which individual analysts mitigate topographical changes on the surface introduces variation to the measurements taken. This can result in major inconsistencies across analysts and even among a single analyst's repeated measurements. Though we did not test inter-observer error, we have demonstrated that the precision and accuracy of the virtual goniometer far surpasses the capabilities of the contact goniometer. Additionally, the virtual goniometer automatically exports the segmentation parameter, a measure of the goodness of fit, the measurement location, radius, and vertices, none of which can be garnered by using the contact goniometer. By providing numerical metadata for the measurement and 3D visualizations, inter- and intra-observer discrepancies can be easily identified and measurements can be replicated with precision. The virtual goniometer resolves the inherent limitations of the contact goniometer and gives researchers a tangible way to discuss how best to employ goniometry to address questions in anthropology.

\section{Conclusion}

The virtual goniometer offers a more precise and consistent way to measure angles and extract data that cannot be extracted manually. Details of the data are recorded such that they can be clearly communicated and replicated. The virtual goniometer provides flexibility that allows the user to adjust parameters and choose an approach that is dependable and useful. Scanning objects is becoming mainstream in anthropology and the virtual goniometer is a logical next step that is easy to integrate into research that uses 3D models. The purpose of this project was to introduce the virtual goniometer and demonstrate its capabilities. The next step for future research is to explore where and how to take measurements on different material types, including analysis of appropriate scan resolution, to address specific anthroplogical questions. \\

\begin{footnotesize}
\begin{noindent}
\textbf{Acknowledgments} 
Thank you to Scott Salonek with the Elk Marketing Council and Christine Kvapil with Crescent Quality Meats for the bones used in this research. We thank the hyenas and their caretakers at the Milwaukee County Zoo and Irvine Park Zoo in Chippewa Falls, Wisconsin and the various math and anthropology student volunteers who broke bones using stone tools. Thank you to Sevin Antley, Chloe Siewart, Mckenzie Sweno, Alexa Krahn, Monica Msechu, Fiona Statz, Emily Sponsel, Kameron Kopps, and Kyra Johnson for helping to break, clean, curate and prepare fragments for scanning. Thank you to Cassandra Koldenhoven and Todd Kes in the Department of Radiology at the Center for Magnetic Resonance Research (CMRR) for CT scanning the fragments. Thank you to Samantha Porter with the University of Minnesota's Advanced Imaging Service for Objects and Spaces (AISOS) who scanned the crystal. Bo Hessburg and Pedro Angulo-Uma\~na worked on the virtual goniometer. Pedro and Carter Chain worked on surfacing the CT scans. Thank you Matt Edling and the University of Minnesota's Evolutionary Anthropology Labs for support in coordinating sessions for bone breakage and guidance for curation. Thank you Abby Brown and the Anatomy Laboratory in the University of Minnesota's College of Veterinary Medicine for providing protocols and a facility to clean bones. Thank you Henry Wyneken and the Liberal Arts Technologies and Innovation Services (LATIS) for statistical consulting.
\end{noindent} \\

\begin{noindent}
\textbf{Funding Information} We would like to thank the National Science Foundation NSF Grant DMS-1816917 and the University of Minnesota's Department of Anthropology for funding this research. 
\end{noindent}
\end{footnotesize} \\


\end{document}